\documentclass[11pt]{article}%
\usepackage{amsfonts}
\usepackage{amssymb}
\usepackage{graphicx}%

\newtheorem{theorem}{Theorem}

\newtheorem{corollary}[theorem]{Corollary}

\newtheorem{lemma}[theorem]{Lemma}

\begin{document}

\begin{center}{\Large An Elementary Proof of the Free-additivity \\
of Voiculescu's Free
Entropy}\end{center}

\vspace{0.5cm}
\begin{center}Don Hadwin \hspace{0.5cm} Weihua Li \hspace{0.5cm}
Junhao Shen
\end{center}
\vspace{0.3cm}

\centerline{Mathematics Department, University of New Hampshire,
Durham, NH 03824}

\vspace{0.3cm}

\centerline{Email:  don@unh.edu \hspace{0.6cm} whli@cisunix.unh.edu
\hspace{0.6cm} jog2@cisunix.unh.edu}

\vspace{0.3cm} \noindent\textbf {Abstract}\ \  D. Voiculescu
\cite{Vo} proved that a standard family of independent random
unitary $k\times k$ matrices and a constant $k\times k$ unitary
matrix is asymtotically free as $k\rightarrow\infty$. This result
was a key ingredient in Voiculescu's proof \cite{Vo1} that his free
entropy is additive when the variables are free. In this paper, we
give a very elementary proof of a more detailed version of this
result \cite{Vo}.

\section{Preliminaries}

\hspace{1.5em}The theory of free probability and free entropy was
introduced by D. Voiculescu in the 1980's, and has become one of the
most powerful and exciting new tools in the theory of von Neumann
algebras. D. Voiculescu \cite{Vo} proved that a standard family of
independent random unitary $k\times k$ matrices and a constant
$k\times k$ unitary matrix is asymtotically free as $k\rightarrow
\infty$. To prove this result, Voiculescu used his noncommutative
central limit theorem and the fact that the unitaries in the polar
decomposition of a family of standard Gaussian random matrices form
a standard family of independent unitary $k\times k$ random
matrices. In this paper, we will give a very elementary proof that
 uses only the basic properties of Haar measure and the definition of a
unitary matrix.

Let $\mathcal{M}_{k}(\mathbb{C)}$ be the $k\times k$ full matrix
algebra with entries in $\mathbb{ C}$, and $\tau_{k}$ be the
normalized trace on $\mathcal{M}_{k}(\mathbb{C)}$, i.e.,
$\tau_{k}=\frac{1}{k}Tr$, where $Tr$ is the usual trace on
$\mathcal{M}_{k}(\mathbb{C)}$. Let ${\mathcal{U}}_{k}$ be the group
of all unitary matrices in $\mathcal{M}_{k}(\mathbb{ C)}$. For
$1\leq i,j\leq k$, define $f_{i,j}: {\mathcal{M}}_{k}({\Bbb
C})\rightarrow\mathbb{C}$ so that any element $a$ in
${\mathcal{M}}_{k}({\Bbb C})$ is the matrix $(f_{i,j}(a));$ i.e.,
$f_{ij}(a)$ is the $(i,j)$-entry of $a$.

An $k\times k$ matrix $u$ is unitary if and only if
\vspace{0.1cm}\newline (1)  $
\sum_{i=1}^{k}|f_{i,j_1}(u)|^{2}=\sum_{j=1}^{k}|f_{i_1,j}(u)|^{2}=1$
for $1\leq i_1, j_1\leq k$, and \vspace{0.1cm}\newline
(2)  $\sum_{i=1}^{k}f_{i,j_{1}}(u)\overline{f_{i,j_{2}}(u)}=\sum_{j=1}^{k}%
f_{i_{1},j}(u)\overline{f_{i_{2},j}(u)}=0,$ whenever $i_{1}\neq
i_{2}$ and $j_{1}\neq j_{2}$.

\vspace{0.1cm} Since ${\mathcal{U}}_{k}$ is a compact group, there
exists a unique normalized Haar measure $\mu_{k}$ on
${\mathcal{U}}_{k}$. In addition,
\[
\int_{{\mathcal{U}}_{k}}f(u)d\mu_{k}(u)=\int_{{\mathcal{U}}_{k}}f(vu)d\mu
_{k}(u)=\int_{{\mathcal{U}}_{k}}f(uv)d\mu_{k}(u).
\]
for every continuous function $f:
{\mathcal{U}}_{k}\rightarrow\mathbb{C}$ and $v\in{\mathcal{U}}_{k}$.

By the tanslation-invariance of $\mu_{k}$, we have the following
lemmas (see also Lemma 12, Lemma 13 and Lemma 14 in \cite{Don}).

\begin{lemma}\label{lemma,permutation}
If $g: \mathbb{C}^{n}\rightarrow\mathbb{C}^{n}$ is a continuous
function, $\sigma$ and $\rho$ are permutations of $\{1,2,\ldots
,k\}$, then
\begin{eqnarray*}
&& \int_{{\mathcal{U}}_{k}}g\left(  f_{i_{1}j_{1}}(u),f_{i_{2}j_{2}
}(u),\ldots, f_{i_{n}j_{n}}(u)\right)  d\mu_{k}(u)\\
&=&\int_{{\mathcal{U}}_{k}}g(f_{\sigma(i_{1}),\rho(j_{1})}(u),f_{\sigma
(i_{2}),\rho(j_{2})}(u),\ldots,
f_{\sigma(i_{n}),\rho(j_{n})}(u))d\mu_{k}(u). \end{eqnarray*}

\end{lemma}

\begin{lemma}\label{lemma, nonzero}
If $\int_{{\mathcal{U}}_{k}} f_{i_{1}j_{1}}(u)\cdots
f_{i_{m}j_{m}}(u) \overline{f_{s_{1}t_{1}}(u)}\cdots\overline{f_{s_{r}%
t_{r}}(u)}d\mu_{k}(u)\neq0$, then

(1) $m=r$,

(2) $(i_{1}, i_{2},\ldots, i_{m})$ is a permutation of $(s_{1},
s_{2},\ldots, s_{r})$,

(3) $(j_{1}, j_{2},\ldots, j_{m})$ is a permutation of $(t_{1},
t_{2},\ldots, t_{r})$.
\end{lemma}

\begin{lemma}\label{lemma,absolute value <}
If $d$ is the maximum cardinality of the sets
$\{i_{1},\ldots,i_{n}\},$ $\{j_{1},\ldots,j_{n}\},$
$\{s_{1},\ldots,s_{r}\}$ and $\{t_{1},\ldots,t_{r}\},$ then, for
every positive integer $k\geq d,$
\[
\left\vert \int_{\mathcal{U}_{k}}f_{i_{1}j_{1}}(u)\cdots
f_{i_{n}j_{n}
}(u)\overline{f_{s_{1}t_{1}}(u)}\cdots\overline{f_{s_{r}t_{r}}(u)}d\mu
_{k}(u)\right\vert \leq\frac{1}{P(k,d)},
\]
where $P\left(  k,d\right)  =k\left(  k-1\right)  \cdots\left(
k-d+1\right) $.
\end{lemma}

\section{Main result}

\hspace{1.5em}If $f: {\cal F}\rightarrow \Bbb C$, let $\left\Vert
f\right\Vert _{\infty}=\sup\left\{ |f(x)|:x\in {\cal F}\right\} .$

\begin{lemma}
\label{0sums}Let $n, m,k$ be positive integers. Let $F,G$ be finite
subsets of $\Bbb N$ with $n=Card(F)$ and $m=Card(G)$. Suppose
$\{f_{i}, g_j: i\in F, j\in G\}$ is a family of mappings from
$\{1,\ldots,k\}$ to $\mathbb{C}$ such that
$\sum_{a=1}^{k}f_{i}\left(a\right) =0$ for $i\in F.$ Let
$H=\{1,\ldots, k\}$. Then
$$ \left\vert \sum_{\sigma:F\cup G\stackrel{1-1}{\rightarrow}H}{\prod _{i\in F}}
 f_{i}\left(  \sigma\left(  i\right)  \right) {\prod
_{j\in G}}  g_{j}\left( \sigma\left( j\right) \right) \right\vert
\leq k^{m+\frac{n}{2}}\left( n+m\right)  ^{n} {\prod _{i\in F}}
\left\Vert f_{i}\right\Vert _{\infty} {\prod _{j\in G}} \left\Vert
g_{j}\right\Vert _{\infty}.
$$
\end{lemma}

Proof.  The proof is by induction on $n.$ When $n=0,$ the obvious
interpretation of the inequality is $$ \left\vert
\sum_{\sigma}{\prod_{j\in G}}
g_{j}\left(  \sigma\left(  j\right)  \right)  \right\vert \leq k^{m}%
{\prod _{j\in G}}  \left\Vert g_{j}\right\Vert _{\infty},
$$
and it holds since the number of functions $\sigma$ is no more than
$k^{m}$.

Suppose the lemma holds for $n$. For $n+1$, let $E$ be a subset of
$F$ with cardinality $n$, say it $E=F\setminus \{b\}$, where $b\in
F$. Then we can define a one-to-one mapping $\sigma:F\cup
G\rightarrow H$ by defining the one-to-one mapping $\sigma: E\cup
G\rightarrow H$ and choosing $s\notin\sigma\left( E\cup G\right) $
to be $\sigma(b)$. Then
\begin{eqnarray*}
&&\left\vert \sum_{\sigma:F\cup G\stackrel{1-1}{\rightarrow}H}{\prod
_{i\in F}}
 f_{i}\left(  \sigma\left(  i\right)  \right) {\prod
_{j\in G}}  g_{j}\left( \sigma\left( j\right) \right)
\right\vert\\
&=&  \left\vert \sum_{\sigma:E\cup
G\stackrel{1-1}{\rightarrow}H}\left( \sum_{s\notin\sigma(E\cup
G)}f_{b}(s)\right) \prod_{i\in E}f_{i}(\sigma
(i))\prod_{j\in G}g_{j}(\sigma(j))\right\vert \\
&=&\left\vert \sum_{\sigma:E\cup
G\stackrel{1-1}{\rightarrow}H}\left( \sum_{s=1}^{k}f_{b}(s)\right)
\prod_{i\in E}f_{i}(\sigma
(i))\prod_{j\in G}g_{j}(\sigma(j))\right. \\
&&  -\left.  \sum_{\sigma:E\cup G\stackrel{1-1}{\rightarrow}H}\left(
\sum _{s\in\sigma(E\cup G)}f_{b}(s)\right) \prod_{i\in
E}f_{i}(\sigma
(i))\prod_{j\in G}g_{j}(\sigma(j))\right\vert \\
&&(\mbox{using}\ \sum_{s=1}^k f_{n+1}(s)=0) \\
&  =&\left\vert \sum_{\sigma:E\cup
G\stackrel{1-1}{\rightarrow}H}\left( \sum_{s\in\sigma(E\cup
G)}f_{b}(s)\right) \prod_{i\in E}f_{i}(\sigma
(i))\prod_{j\in G}g_{j}(\sigma(j))\right\vert \\
&  =&\left\vert \sum_{\sigma:E\cup
G\stackrel{1-1}{\rightarrow}H}\left( \sum_{t\in E\cup
G}f_{b}(\sigma(t))\right) \prod_{i\in E}f_{i}(\sigma
(i))\prod_{j\in G}g_{j}(\sigma(j))\right\vert \\
 &  \leq&\left|
\sum_{\sigma:E\cup G\stackrel{1-1}{\rightarrow}H}\left( \sum_{t\in
E}f_{b}(\sigma(t))\right) \prod_{i\in E}f_{i}(\sigma
(i))\prod_{j\in G}g_{j}(\sigma(j))\right| \\
&&  +\left\vert \sum_{\sigma:E\cup
G\stackrel{1-1}{\rightarrow}H}\left( \sum_{t\in
G}f_{b}(\sigma(t))\right) \prod_{i\in E}f_{i}
(\sigma(i))\prod_{j\in G}g_{j}(\sigma(j))\right\vert \\
&  \leq &\sum_{t\in E}\left\vert \sum_{\sigma:E\cup
G\stackrel{1-1}{\rightarrow} H}\left( \prod_{i\in E,i\neq
t}f_{i}(\sigma(i))\right)  (f_{b}
f_{t})(\sigma(t))\prod_{j\in G}g_{j}(\sigma(j))\right\vert \\
&&  +\sum_{t\in G}\left\vert \sum_{\sigma:E\cup G\stackrel{1-1}
{\rightarrow}H}\left(  \prod_{i\in E,i\neq t}f_{i}(\sigma(i))\right)
(f_{b}f_{t})(\sigma(t))\prod_{j\in G}g_{j}(\sigma
(j))\right\vert \\
&&(\mbox{using induction on the quantities inside the absolute value signs}\\
&&\mbox{ and viewing }\ f_{b}f_t \ \mbox{as a single function})\\
 &  \leq&
n((m+1)+(n-1))^{n-1}k^{\frac{n-1}{2}+m+1}\prod_{i\in F}\Vert
f_{i}\Vert_{\infty}\prod_{j\in G}\Vert g_{j}\Vert_{\infty}\\
&&  +m(m+n)^{n}k^{\frac{n}{2}+m}\prod_{i\in F}\Vert
f_{i}\Vert_{\infty}
\prod_{j\in G}\Vert g_{j}\Vert_{\infty}\\
&  \leq&(m+n+1)^{n+1}k^{\frac{n+1}{2}+m}\prod_{i\in F}\Vert f_{i}
\Vert_{\infty}\prod_{j\in G}\Vert g_{j}\Vert_{\infty}.
\end{eqnarray*}\hfill$\Box$

\vspace{0.2cm}
 Let $\mathcal{U}_{k}^{n}$ denote the direct product of $n$ copies
of $\mathcal{U}_{k}$, and $\mu_{k}^{n}$ denote the corresponding
product measure. Let $C\left(
\mathcal{U}_{k}^{n},\mathcal{M}_{k}(\mathbb{C)}\right)  $ denote the
C$^{\ast}$-algebra of all continuous functions from
$\mathcal{U}_{k}^{n}$ into $\mathcal{M}_{k}(\mathbb{C)}$. If
$\overrightarrow{u}=\left( u_{1},\ldots,u_{n}\right)
\in\mathcal{U}_{k}^{n}$, then the coordinate variables
$u_{1},\ldots,u_{n}$ are unitary elements of $C\left(
\mathcal{U}_{k}^{n},\mathcal{M}_{k}(\mathbb{C)}\right) $.

The following lemma is a vastly improved estimate that is
independent of the maximum cardinality of the indices in the
integral. We require the elementary inequalities
$m^{m}\leq2^{m^{2}}$ and $\frac{1}{P\left(k,m\right)}
\leq\frac{m^{m}}{k^{m}}$ for positive integers $m\leq k.$

\begin{lemma}
\label{lemma,nuest}Suppose $m$ is a positive integer. For every
positive integers $k,n$ with $k\geq m,$ and for all subsets
$\{i_{1},\ldots,i_{m}\}$, $\{j_{1},\ldots,j_{m}\}$ of $\left\{
1,\ldots,k\right\}$, and $\{\iota_{1},\ldots,\iota_{m},
\eta_{1}\ldots,\eta_{m}\}$ of $\{ 1,\ldots,n\}$,
\[
\left\vert
\int_{\mathcal{U}_{k}^{n}}f_{i_{1}j_{1}}(u_{\iota_{1}})\cdots
f_{i_{m}j_{m}}(u_{\iota_{m}})\overline{f_{s_{1}t_{1}}(u_{\eta_{1}})}\cdots
\overline{f_{s_{m}t_{m}}(u_{\eta_{m}})}d\mu_{k}^{n}(\vec{u})\right\vert
\leq\frac{4^{m^{2}}}{k^{m}}.
\]

\end{lemma}

Proof. for $1\leq j\leq n$, let $T_j=\{1\leq \lambda\leq m:
\iota_{\lambda}=j\}$ and $T_j^*=\{1\leq \lambda\leq m:
\eta_{\lambda}=j\}$. Then \begin{eqnarray*} &&
\int_{\mathcal{U}_{k}^{n}}f_{i_{1}j_{1}}(u_{\iota_{1}})\cdots
f_{i_{m}j_{m}}(u_{\iota_{m}})\overline{f_{s_{1}t_{1}}(u_{\eta_{1}})}\cdots
\overline{f_{s_{m}t_{m}}(u_{\eta_{m}})}d\mu_{k}^{n}(\vec{u})\\
&=&\prod_{j=1}^n\int_{{\cal U}_k}(\prod_{\lambda\in
T_j}f_{i_{\lambda}j_{\lambda}}(u_j)\prod_{\lambda\in
T_j^*}\overline{f_{s_{\lambda}
t_{\lambda}}(u_j)})d\mu_k(u_j).\end{eqnarray*} Hence, we can assume
that $n=1$. Moreover, in view of the Cauchy-Schwarz inequality, it
is sufficient to prove that
\begin{eqnarray} I=\int_{{\cal U}_{k}}\left\vert f_{i_{1}j_{1}}\left(
u\right) \right\vert ^{2}\left\vert f_{i_{2}j_{2}}\left(  u\right)
\right\vert ^{2}\cdots \left\vert f_{i_{m}j_{m}}\left(  u\right)
\right\vert ^{2}d\mu_{k}\left( u\right)\leq  \frac{4^{m^{2}}}{k^{m}}
.
\end{eqnarray}

Let $d$ be the maximum cardinality of the sets
$\{i_{1},\ldots,i_{m}\}$ and $\{j_{1},\ldots,j_{m}\}$. By replacing
$u$ with $u^{\ast},$ which does not alter the integral but
interchanges $i$'s with $j$'s, we can assume that $d$ is the
cardinality of $\{i_{1},\ldots,i_{m}\}$ . Then $1\leq d\leq m.$ Let
$B_{d,k}$ be the largest integral of the type in (1) with
$d=Card(\{i_{1},\ldots,i_{m}\}).$

If $d=m,$ then, by Lemma \ref{lemma,absolute value <}, the integral
in (1) is at most $\frac{1}{P(k,m)}$, and $\frac{1}{P(k,m)}\leq
\frac{m^m}{k^m}\leq \frac{4^{m^2}}{k^{m}}.$

Now we will prove that $B_{d,k}\leq 2^mB_{d+1,k}$ whenever $1\leq d<
m.$ For $1\leq d<m$, assume that the integral $I$ above equals
$B_{d,k}$. Since $d<m,$ at least two of $i_{1},\ldots,i_{m}$ must be
the same. From Lemma \ref{lemma,permutation}, we can assume that
$1\leq i_{1},\ldots,i_{m}\leq d$ and $1=i_{1}=i_{2}=\cdots=i_{s}$
and $1\notin\left\{ i_{s+1},\ldots,i_{m}\right\}  .$ Since $k\geq
m>d,$ we can define a unitary matrix $v$ with $1$ on the diagonal
except in the $\left( 1,1\right)  $ and $\left(  k,k\right)  $
positions, with $\frac{1}{\sqrt{2}}$ in the $\left(  1,1\right)
,\left( k,1\right) ,\left(  k,k\right)  $ positions and
$-\frac{1}{\sqrt{2}}$ in the $\left(  1,k\right)  $ position. Since
the integral remains unchanged when we replace the variable $u$ with
$vu,$ we obtain
\begin{eqnarray*}
B_{d,k}&=&\frac{1}{2^{s}}\int_{{\cal U}_{k}} {\prod _{\beta=1}^{s}}
\left\vert f_{1j_{\beta}}\left(  u\right)
+f_{kj_{\beta}}(u)\right\vert ^{2} {\prod _{\alpha=s+1}^{m}}
  \left\vert f_{i_{\alpha}j_{\alpha}}\left(  u\right) \right\vert
^{2}d\mu _{k}\left(  u\right)\\
&=&\frac{1}{2^{s}}\int_{{\cal U}_{k}} {\prod _{\beta=1}^{s}}
 \left(  \left\vert f_{1j_{\beta}}\left(  u\right) \right\vert
^{2}+\overline {f_{1j_{\beta}}\left(  u\right)
}f_{kj_{\beta}}(u)+f_{1j_{\beta}}\left(  u\right)
\overline{f_{kj_{\beta}}(u)}+\left\vert f_{kj_{\beta}}\left(
u\right)
\right\vert ^{2}\right)\cdot\\
&&{\prod _{\alpha=s+1}^{m}}
 \left\vert f_{i_{\alpha}j_{\alpha}}\left(  u\right) \right\vert
^{2}d\mu
_{k}\left(  u\right) \\
& =&\frac{1}{2^{s}}\int_{{\cal U}_{k}} {\prod_{\beta=1}^{s}}
\left\vert f_{1j_{\beta}}\left(  u\right)  \right\vert ^{2} {\prod
_{\alpha=s+1}^{m}}  \left\vert f_{i_{\alpha}j_{\alpha}}\left(
u\right) \right\vert ^{2}d\mu
_{k}\left(  u\right) \\
&&  +\frac{1}{2^{s}}\int_{{\cal U}_{k}} {\prod _{\beta=1}^{s}}
\left\vert f_{kj_{\beta}}\left(  u\right) \right\vert ^{2} {\prod
_{\alpha=s+1}^{m}}
 \left\vert f_{i_{\alpha}j_{\alpha}}\left(  u\right) \right\vert
^{2}d\mu _{k}\left(  u\right) +\frac{1}{2^{s}}\int_{{\cal
U}_{k}}\Delta d\mu_{k}\left( u\right),
\end{eqnarray*}
where $\Delta$ is a summation of $4^{s}-2$ terms with each of them
having both an $f_{1\ast}\left(  u\right)  $ and an $f_{k\ast}\left(
u\right)  $ factor (with or without conjugation signs) and the
maximum cardinality of the indices in each term is $d+1$, which
implies $\left|\int_{{\cal U}_{k}}\Delta d\mu_{k}\left(
u\right)\right|\leq (4^s-2)B_{d+1,k}$.

Since \begin{eqnarray*} B_{d,k}&=&\int_{{\cal U}_{k}}
{\prod_{\beta=1}^{s}} \left\vert f_{1j_{\beta}}\left(  u\right)
\right\vert ^{2} {\prod_{\alpha=s+1}^{m}} \left\vert
f_{i_{\alpha}j_{\alpha}}\left(
u\right)  \right\vert ^{2}d\mu _{k}\left(  u\right)\\
&=&\int_{{\cal U}_{k}} {\prod_{\beta=1}^{s}}
\left\vert f_{kj_{\beta}}\left(  u\right) \right\vert ^{2}%
{\prod_{\alpha=s+1}^{m}} \left\vert f_{i_{\alpha}j_{\alpha}}\left(
u\right)  \right\vert ^{2}d\mu _{k}\left(u\right),
\end{eqnarray*}
we have $$B_{d,k}\leq \frac{1}{2^{s}}\left(
B_{d,k}+B_{d,k}\right)+\frac{1} {2^{s}}\left( 4^{s}-2\right)
B_{d+1,k}.$$ Therefore
$$
B_{d,k}\leq 2^{m}B_{d+1,k}.
$$
It follows that $B_{d,k}\leq2^{m\left(  m-d\right) }B_{m,k}\leq\frac
{2^{m^{2}}}{P\left(k,m\right) }\leq\frac{2^{m^{2}}m^{m}}{k^{m}}\leq
\frac{4^{m^{2}}}{k^{m}}$ when $k\geq m$ and $1\leq d\leq m.$
\hfill$\Box$

\vspace{0.3cm}

For any positive integer $m$, let $B\left(  m\right)  $ be the
\emph{Bell number} of $m,$ i.e., the number of equivalence relations
on a set with cardinality $m$. Suppose $\cal M$ is a von Neumann
algebra with a faithful tracial state $\tau$ and ${\cal U(M)}$ is
the set of all unitary elements in $\cal M$ and
$\overrightarrow{u}=(u_1,\ldots,u_n)\in {\cal U}({\cal M})^n$. Let
${\Bbb F}_n$ be a free group with standard generators $h_1, \ldots,
h_n$. Then there is a homomorphism $\rho:$
$\mathbb{F}_{n}\rightarrow\mathcal{U(M)}$ such that
$\rho(h_{j})=u_{j}.$ We use the notation $\rho(g)=g(\overrightarrow
{u})=g(u_{1},\ldots,u_{n}).$

D. Voiculescu \cite{Vo} proved that a standard family of independent
random unitary $k\times k$ matrices and a constant $k\times k$
unitary matrix is asymtotically free as $k\rightarrow\infty$. The
following theorem gives a very elementary proof of a more detailed
version of D. Voiculescu's result. The constants in the following
theorem are far from best possible, but they are, at least,
explicit.

\begin{theorem}\label{theorem,estimation}
Suppose $M>0$ and $m$ is a positive integer. For every reduced words
$g_{1},\ldots,g_{w}\in\mathbb{F}_{n}\backslash\{e\}$ with
$\sum_{i=1}^{w} length\left( g_{i}\right) =m$, and commuting normal
$k\times k$ matrices $x_{1},\ldots,x_{w}$ with trace 0 and
$\left\Vert x_{i}\right\Vert \leq M$ for all $1\leq i\leq w$, we
have

\begin{enumerate}
\item
\[
\left\vert \int_{\mathcal{U}_{k}^{n}}\tau_{k}\left(  g_{1}\left(
\vec {u}\right)  x_{1}g_{2}\left(  \vec{u}\right)  x_{2}\cdots
g_{w}\left(  \vec {u}\right)  x_{w}\right)  d\mu_{k}^{n}\left(
\vec{u}\right)  \right\vert
\leq\frac{B\left(  m\right)  2^{m^{2}}\left(  Mw\right)  ^{w}}{k},%
\]

\item
\[
\int_{\mathcal{U}_{k}^{n}}\left\vert \tau_{k}\left(  g_{1}\left(
\vec {u}\right)  x_{1}g_{2}\left(  \vec{u}\right)  x_{2}\cdots
g_{w}\left(  \vec {u}\right)  x_{w}\right)  \right\vert
^{2}d\mu_{k}^{n}\left(  \vec{u}\right)
\leq\frac{B(2m)4^{m^2}(2Mw)^{2w}}{k^{2}},%
\]

\item
if $\varepsilon>0$, and $k>\frac{2B\left(  m\right) 2^{m^{2}}\left(
Mw\right)  ^{w}}{\varepsilon}$, then
\begin{eqnarray*}\mu_k^n\left(\left\{\overrightarrow{v}\in {\cal U}_k^n:
\left|\tau_{k}\left(  g_{1}\left( \vec {v}\right) x_{1}g_{2}\left(
\vec{v}\right)  x_{2}\cdots g_{w}\left( \vec {v}\right)
x_{w}\right)\right|\geq \varepsilon\right\}\right)\leq
\frac{4B(2m)4^{m^2}(2M\omega)^{2\omega}}{k^2\varepsilon^2}.\end{eqnarray*}
\end{enumerate}
\end{theorem}

Proof. Since $x_{1},\ldots,x_{w}$ are commuting normal matrices,
there is a unitary matrix $v$ such that, for $1\leq j\leq w$,
$vx_jv^*=a_j$ is diagonal. Since $\tau_k$ is tracial and
$$g_{1}\left( \vec {u}\right) x_{1}g_{2}\left( \vec{u}\right)
x_{2}\cdots g_{w}\left(  \vec {u}\right) x_{w}=v^*\left( g_{1}\left(
v\vec {u}v^*\right) a_1g_{2}\left( v\vec{u}v^*\right) a_2\cdots
g_{w}\left( v\vec {u}v^*\right) a_w\right)v,$$ we have
$$\tau_{k}\left( g_{1}\left( \vec {u}\right) x_{1}g_{2}\left(
\vec{u}\right) x_{2}\cdots g_{w}\left(  \vec {u}\right)
x_{w}\right)=\tau_{k}\left( g_{1}\left( v\vec {u}v^*\right)
a_1g_{2}\left( v\vec{u}v^*\right) a_2\cdots g_{w}\left( v\vec
{u}v^*\right) a_w\right).
$$
By the translation-invariance of $\mu_k^n$, we can assume that $x_1,
\ldots, x_w$ are all diagonal matrices.

{\it Proof of the first statement}. Write $g_{1}( \vec
{u})=u_{s_{1}}^{\varepsilon_{1}}\cdots
u_{s_{m_{1}}}^{\varepsilon_{m_{1}}},$ $g_{2}( \vec
{u})=u_{s_{m_{1}+1}}^{\varepsilon _{m_{1}+1}}\cdots
u_{s_{m_{2}}}^{\varepsilon_{m_{2}}},\ldots,$ $g_{w}( \vec {u})
=u_{s_{m_{w-1}+1}}^{\varepsilon_{m_{w-1}+1}}\cdots u_{s_{m_{w}}}%
^{\varepsilon_{m_{w}}}$ with each $\varepsilon_{j}\in\left\{
-1,1\right\}$ and $s_{j}\in\{1,\ldots,n\}$ and with the property
that $s_{j}=s_{j+1}$
implies $\varepsilon_{j}=\varepsilon_{j+1}$ unless $j\in\left\{  m_{1}%
,\ldots,m_{w}\right\}  $. Note that $m_{w}=m$ since $\sum
length\left( g_{i}\right)  =m.$ Also write $x_{j}=diag\left(
\gamma_{j}\left(  1\right) ,\ldots,\gamma_{j}\left(  k\right)
\right) $ for $1\leq j\leq w.$

Define $\dot{+}$ on $\{1, \ldots, m_w=m\}$ by
$s\dot{+}1=\left\{\begin{array}{ll}1,&
s=m_w\\
s+1, & 1\leq s\leq m_w-1\end{array}\right..$ Then we
have\begin{eqnarray*} &&\int_{\mathcal{U}_{k}^{n}}\tau_{k}\left(
g_{1}\left( \vec{u}\right) x_{1}g_{2}\left(  \vec{u}\right)
x_{2}\cdots g_{w}\left( \vec{u}\right) x_{w}\right)
d\mu_{k}^{n}\left(
\vec{u}\right) \\
& =& \frac{1}{k}\sum_{1\leq i_{1},\ldots,i_{m_{w}\dot{+}1}=i_{1}\leq
k}\left( {\prod_{\nu=1}^{w}}
\gamma_{\nu}\left(  i_{m_{\nu}\dot{+}1}\right)  \right)  \int_{\mathcal{U}_{k}^{n}}%
{\prod_{j=1}^{m_{w}}}
f_{i_{j}i_{j\dot{+}1}}\left(  u_{s_{j}}^{\varepsilon_{j}}\right)  d\mu_{k}%
^{n}\left(  \vec{u}\right)  .
\end{eqnarray*}

Let $E=\left\{1,2,\ldots,m_{w}\right\}$. We can represent a choice
of $1\leq i_{1},\ldots,i_{m_{w}}\leq k$ by a function
$\alpha:E\rightarrow H=\{1,\ldots,k\}.$ Thus we can replace the sum
$\sum\limits_{1\leq i_{1},\ldots,i_{m_w\dot{+}1}=i_{1}\leq k}$ with
$\sum\limits_{\alpha:E\rightarrow H}$ in the above equation. That is
$$\left(I=\right)\frac{1}{k}\sum_{\alpha:E\rightarrow H}\left( {\prod_{\nu=1}^{w}}
\gamma_{\nu}\left(  \alpha(m_{\nu}\dot{+}1)\right)  \right)
\int_{\mathcal{U}_{k}^{n}} {\prod_{j=1}^{m_{w}}}
f_{\alpha(j),\alpha(j\dot{+}1)}\left(
u_{s_{j}}^{\varepsilon_{j}}\right) d\mu_{k}^{n}\left(
\vec{u}\right).$$

We only need to restrict sums to the functions $\alpha$ such that
the integral $$ I\left(\alpha\right) =\int_{\mathcal{U}_{k}^{n}}
{\prod_{j=1}^{m_{w}}}
f_{\alpha\left(j\right),\alpha\left(j\dot{+}1\right)}\left(u_{s_{j}}^{\varepsilon_{j}}\right)
d\mu_{k}^{n}\left(\vec{u}\right) \neq0.$$ We call such function
$\alpha$ $\emph{good}$, thus $$I=
\frac{1}{k}\sum_{\begin{array}{l}\alpha:E\rightarrow H\\
\alpha\ \mbox{is good}\end{array}}\left( {\prod_{\nu=1}^{w}}
\gamma_{\nu}\left( \alpha(m_{\nu}\dot{+}1)\right)  \right)
I(\alpha).$$Lemma \ref{lemma, nonzero} tells us that $m_{w}$ must be
even and exactly half of the $\varepsilon_{j}$'s are $1$ and the
other half are $-1.$  We know from Lemma \ref{lemma,nuest} that
\begin{equation} \left\vert I\left( \alpha\right) \right\vert
\leq\frac{4^{\left( m/2\right)
^{2}}}{k^{m/2}}\leq\frac{2^{m^{2}}}{k^{m/2}}.
\end{equation} Moreover, Lemma \ref{lemma, nonzero} says that
if $j\in E$ but $j\notin\left\{
1\dot{+}m_{1},\ldots,1\dot{+}m_{w}\right\} ,$ then $\alpha\left(
j\right) =\alpha\left(  j^{\prime}\right)  $ for some
$j^{\prime}\neq j.$

Next we define an equivalence relation $\sim_{\alpha}$ on $E$ by
saying $i\sim_{\alpha}j$ if and only if $\alpha\left( i\right)
=\alpha\left( j\right)  .$ Note that if $\beta:E\rightarrow H,$ then
the relations $\sim_{\alpha}$ and $\sim_{\beta}$ are equal if and
only if there is a permutation $\sigma:H\rightarrow H$ such that
$\beta=\sigma\circ\alpha.$ We define an equivalent relation
$\approx$ on the set of all good functions by
$$\alpha\approx\beta \ \mbox{if and only if}\ \sim_{\alpha}=\sim_{\beta}.$$ It is clear that $$\alpha\approx\beta\Longrightarrow I\left( \alpha\right)
=I\left(  \beta\right).$$

 If
$j\in E,$ let $\left[  j\right]  _{\alpha}$ denote the $\sim_{\alpha}%
$-equivalence class of $j$, and let $E_{\alpha}$ denote the set of
all such equivalence classes. We can construct all of the functions
$\beta$ equivalent
to $\alpha$ in terms of injective functions%
\[
\sigma:E_{\alpha}\stackrel{1\mbox{-}1}{\rightarrow}H
\]
by defining
\[
\beta\left(  j\right)  =\sigma\left(  \left[  j\right]
_{\alpha}\right)  .
\]

Let $A$ contains exactly one function $\alpha$ from each
$\approx$-equivalence class of good functions. Then we can
write\begin{eqnarray}
\left|I\right|&=&\left|\frac{1}{k}\sum_{\begin{array}{l}\alpha:E\rightarrow
H\\ \alpha\ \mbox{is good}\end{array}}\left( {\prod_{\nu=1}^{w}}
\gamma_{\nu}\left(
\alpha(m_{\nu}\dot{+}1)\right)  \right) I(\alpha)\right|\nonumber\\
&=& \left|\frac{1}{k}\sum_{\alpha\in A}
I(\alpha)\sum_{\beta\approx\alpha}{\prod_{\nu=1}^{w}}
\gamma_{\nu}\left( \beta\left(  m_{\nu}\dot{+}1
\right) \right)\right|\nonumber\\
&=&\frac{1}{k}\left|\sum_{\alpha\in A}|I\left(
\alpha\right)\sum_{\sigma:E_{\alpha
}\stackrel{1\mbox{-}1}{\rightarrow}H} {\prod_{\nu=1}^{w}}
\gamma_{\nu}\left( \sigma\left( \left[ m_{\nu}\dot{+}1\right]
_{\alpha}\right) \right)\right|\nonumber\\
&\leq&\frac{1}{k}\sum_{\alpha\in A}\left|I\left(
\alpha\right)\right|\left|\sum_{\sigma:E_{\alpha
}\stackrel{1\mbox{-}1}{\rightarrow}H} {\prod_{\nu=1}^{w}}
\gamma_{\nu}\left( \sigma\left( \left[ m_ {\nu}\dot{+}1\right]
_{\alpha}\right) \right)\right|.\end{eqnarray}

Also we know that \begin{equation}Card\left(  A\right)  \leq B\left(
m\right) .\end{equation}

We only need to focus on $
\left|\sum_{\sigma:E_{\alpha}\stackrel{1\mbox{-}1}{\rightarrow}H}%
{\prod_{\nu=1}^{w}}  \gamma_{\nu}\left(  \sigma\left( \left[
 m_ {\nu}\dot{+}1\right]  _{\alpha}\right) \right)\right|. $ Let
$$F_{\alpha}=\left\{\left[   m_ {\nu}\dot{+}1\right] _{\alpha}:1\leq\nu\leq
w,Card\left(\left[ m_ {\nu}\dot{+}1\right] _{\alpha}\right)
=1\right\},$$
$$G_{\alpha}=\left\{ \left[   m_ {\nu}\dot{+}1\right] _{\alpha}:1\leq
\nu\leq w,Card\left( \left[   m_ {\nu}\dot{+}1\right]
_{\alpha}\right)
>1\right\},$$  $$K_{\alpha}=E_{\alpha}\backslash\left(
F_{\alpha}\cup G_{\alpha }\right).$$ Since the product
${\prod_{\nu=1}^{w}} \gamma_{\nu}\left(  \sigma\left(  \left[
 m_ {\nu}\dot{+}1\right] _{\alpha}\right) \right) $ is determined once
$\sigma$ is defined on $F_{\alpha}\cup G_{\alpha},$ it follows that
this product is repeated at most $P\left(k,card\left(
K_{\alpha}\right)  \right)  $ times. Hence we
have\begin{eqnarray}&&\left|\sum_{\sigma:E_{\alpha}\stackrel{1\mbox{-}1}{\rightarrow}H}
{\prod_{\nu=1}^{w}}  \gamma_{\nu}\left(  \sigma\left( \left[
 m_ {\nu}\dot{+}1\right]  _{\alpha}\right) \right)\right| \nonumber\\
&\leq & P\left(k,card\left( K_{\alpha}\right) \right)\left| \sum
_{\sigma:F_{\alpha}\cup
G_{\alpha}\stackrel{1\mbox{-}1}{\rightarrow}F}{\prod_{\nu=1}^{w}}
\gamma_{\nu}\left(  \sigma\left( \left[  m_ {\nu}\dot{+}1\right]
_{\alpha}\right) \right)\right|\nonumber\\
 &\leq&k^{card\left( K_{\alpha}\right)  }\left| \sum
_{\sigma:F_{\alpha}\cup
G_{\alpha}\stackrel{1\mbox{-}1}{\rightarrow}F}{\prod_{\nu=1}^{w}}
\gamma_{\nu}\left(  \sigma\left( \left[  m_ {\nu}\dot{+}1\right]
_{\alpha}\right) \right)\right|.
\end{eqnarray}

If $a=[ m_ {\nu}\dot{+}1]_{\alpha}\in F_{\alpha}$, from the
definition of $F_{\alpha}$, it is clear that $\nu$ is unique. Then
define $f_a\left(  \sigma\left( a\right)\right)=\gamma_{\nu}\left(
\sigma\left( a\right)\right)$. By $\tau_k(x_i)=0$ for all $1\leq
i\leq w$, it follows that $\sum_{s=1}^k f_a(s)=0$. If $b=[ m_
{\nu}\dot{+}1]_{\alpha}\in G_{\alpha}$, from the definition of
$G_{\alpha}$, the cardinality of $b$ is greater than 1, say it $r$.
Then define $g_b\left(
\sigma\left(b\right)\right)=\left(\gamma_{\nu} \left(\sigma\left(
b\right)\right)\right)^r$. Therefore
\begin{eqnarray}&&\left|\sum _{\sigma:F_{\alpha}\cup
G_{\alpha}\stackrel{1\mbox{-}1}{\rightarrow}H}{\prod_{\nu=1}^{w}}
\gamma_{\nu}\left(  \sigma\left( \left[  m_ {\nu}\dot{+}1\right]
_{\alpha}\right) \right)\right|\nonumber\\
&=&\left|\sum _{\sigma:F_{\alpha}\cup
G_{\alpha}\stackrel{1\mbox{-}1}{\rightarrow}H}{\prod_{a\in
F_{\alpha}}}f_a\left(  \sigma\left( a\right)\right){\prod_{b\in
G_{\alpha}}}g_b\left( \sigma\left(
b\right)\right) \right|\nonumber\\
&&\left(\mbox{letting}\ F=F_{\alpha}, G=G_{\alpha}\ \mbox{and using Lemma \ref{0sums} }\right)\nonumber\\
&\leq & k^{\left[  card\left( F_{\alpha}\right)  /2\right]
+card\left( G_{\alpha}\right) }w^{w}M^{w}.
\end{eqnarray}
As we mentioned before that $card\left( \left[ j\right]
_{\alpha}\right) =1$ implies $\left[  j\right] _{\alpha}\in
F_{\alpha},$ we see that
\begin{equation} \left[  card\left( F_{\alpha}\right) /2\right] +card\left(
G_{\alpha }\right) +card\left( K_{\alpha}\right)  \leq card\left(
E\right) /2=m_{w}/2.
\end{equation}
Combining (2), (3), (4), (5), (6) and (7) together, we have $$
|I|\leq \frac{1}{k}B\left(  m\right)  2^{m^{2}}\left( Mw\right)
^{w}.$$

 {\it Proof of the second statement.} We know that
\begin{eqnarray*}&&\left|\tau_{k}\left( g_{1}\left(
\vec{u}\right) x_{1}g_{2}\left(  \vec{u}\right) x_{2}\cdots
g_{w}\left( \vec{u}\right)
x_{w}\right)\right|^2\\
&=&\tau_{k}\left( g_{1}\left( \vec{u}\right) x_{1}g_{2}\left(
\vec{u}\right) x_{2}\cdots g_{w}\left( \vec{u}\right)
x_{w}\right)\cdot \overline{\tau_{k}\left( g_{1}\left(
\vec{u}\right) x_{1}g_{2}\left(  \vec{u}\right) x_{2}\cdots
g_{w}\left( \vec{u}\right) x_{w}\right)}\\
&=&\frac{1}{k^2}\sum_{1\leq i_1, \ldots, i_{m_w+1}=i_1\leq k}\left(
{\prod_{\nu=1}^{w}} \gamma_{\nu}\left(  i_{m_{\nu}+1}\right) \right)
{\prod_{j=1}^{m_{w}}} f_{i_{j}i_{j+1}}\left(
u_{s_{j}}^{\varepsilon_{j}}\right)  \\
&&\sum_{1\leq l_1, \ldots, l_{m_w+1}=l_1\leq k}\left(
{\prod_{\lambda=1}^{w}}\overline{ \gamma_{\lambda}\left(
l_{m_{\lambda}+1}\right) } \right) {\prod_{t=1}^{m_{w}}}
\overline{f_{l_{t}l_{t+1}}\left(
u_{s_{t}}^{\varepsilon_{t}}\right)}.
\end{eqnarray*}

Define $\dot{+}$ on the set $\{1,2,\ldots, 2m_w\}$ by
$x\dot{+}1=\left\{\begin{array}{ll}1,&
x=m_w\\
m_w+1,& x=2m_w\\
 x+1, & 1\leq x\leq m_w-1\end{array}\right..$ Then
we have\begin{eqnarray*}
I&=&\int_{\mathcal{U}_{k}^{n}}\left|\tau_{k}\left( g_{1}\left(
\vec{u}\right) x_{1}g_{2}\left(  \vec{u}\right) x_{2}\cdots
g_{w}\left( \vec{u}\right) x_{w}\right)\right|^2 d\mu_{k}^{n}\left(
\vec{u}\right) \\
&=&\frac{1}{k^2}\int_{\mathcal{U}_{k}^{n}}\sum_{1\leq i_1, \ldots,
i_{m_w\dot{+}1}\leq k}\left( {\prod_{\nu=1}^{w}} \gamma_{\nu}\left(
i_{m_{\nu}\dot{+}1}\right) \right) {\prod_{j=1}^{m_{w}}}
f_{i_{j}i_{j\dot{+}1}}\left(
u_{s_{j}}^{\varepsilon_{j}}\right)  \\
&&\sum_{1\leq l_1, \ldots, l_{m_w\dot{+}1}\leq k}\left(
{\prod_{\lambda=1}^{w}}\overline{ \gamma_{\lambda}\left(
l_{m_{\lambda}\dot{+}1}\right) } \right) {\prod_{t=1}^{m_{w}}}
\overline{f_{l_{t}l_{t\dot{+}1}}\left(
u_{s_{t}}^{\varepsilon_{t}}\right)}.
\end{eqnarray*}

Let $E=\left\{1,2,\ldots,2m_{w}\right\}$. We can represent a choice
of $1\leq i_{1},\ldots,i_{m_{w}}\leq k$ by a function
$\alpha:E\rightarrow H=\{1,\ldots,k\}.$ Thus we can rewrite $I$
\begin{eqnarray*}I&=& \frac{1}{k^2}\sum_{\alpha: E\rightarrow H}\left( {\prod_{\nu=1}^{w}} \gamma_{\nu}\left(
\alpha({m_{\nu}\dot{+}1})\right) \right) \left(
{\prod_{\lambda=1}^{w}}\overline{ \gamma_{\lambda}\left(
\alpha((m_{\lambda}\dot{+}1)+m_w)\right) }
\right)\\
&& \int_{\mathcal{U}_{k}^{n}}{\prod_{j=1}^{m_{w}}}
f_{\alpha(j)\alpha(j\dot{+}1)}\left(
u_{s_{j}}^{\varepsilon_{j}}\right)  {\prod_{t=1}^{m_{w}}}
\overline{f_{\alpha(t+m_w)\alpha((t\dot{+}1)+m_w)}\left(
u_{s_{t}}^{\varepsilon_{t}}\right)}.\end{eqnarray*}

We only need to restrict sums to the functions $\alpha$ such that
the integral $$ I\left(\alpha\right) =
\int_{\mathcal{U}_{k}^{n}}{\prod_{j=1}^{m_{w}}}
f_{\alpha(j)\alpha(j\dot{+}1)}\left(
u_{s_{j}}^{\varepsilon_{j}}\right)  {\prod_{t=1}^{m_{w}}}
\overline{f_{\alpha(t+m_w)\alpha((t\dot{+}1)+m_w)}\left(
u_{s_{t}}^{\varepsilon_{t}}\right)} \neq0.$$ We call such function
$\alpha$ $\emph{good}$. We have \begin{equation}I=
\frac{1}{k^2}\sum_{\alpha: E\rightarrow H}\left( {\prod_{\nu=1}^{w}}
\gamma_{\nu}\left( \alpha({m_{\nu}\dot{+}1})\right) \right) \left(
{\prod_{\lambda=1}^{w}}\overline{ \gamma_{\lambda}\left(
\alpha((m_{\lambda}\dot{+}1)+m_w)\right) }
\right)I(\alpha).\end{equation} Lemma \ref{lemma, nonzero} tells us
that if $j\in E$ but $j\notin\left\{
1\dot{+}m_{1},\ldots,1\dot{+}m_{w}, (1\dot{+}m_1)+m_w, \ldots,
(1\dot{+}m_{w})+m_w\right\} ,$ then $\alpha\left( j\right)
=\alpha\left( j^{\prime}\right)  $ for some $j^{\prime}\neq j.$  We
know from Lemma \ref{lemma,nuest} that \begin{equation} \left\vert
I\left( \alpha\right) \right\vert \leq\frac{4^{\left( m\right)
^{2}}}{k^{m}}.
\end{equation}

Next we define an equivalence relation $\sim_{\alpha}$ on $E$ by
saying $i\sim_{\alpha}j$ if and only if $\alpha\left( i\right)
=\alpha\left( j\right)  .$ Note that if $\beta:E\rightarrow H,$ then
the relations $\sim_{\alpha}$ and $\sim_{\beta}$ are equal if and
only if there is a permutation $\sigma:H\rightarrow H$ such that
$\beta=\sigma\circ\alpha.$ We define an equivalent relation
$\approx$ on the set of all good functions by
$$\alpha\approx\beta \ \mbox{if and only if}\ \sim_{\alpha}=\sim_{\beta}.$$ It is clear that $$\alpha\approx\beta\Longrightarrow I\left( \alpha\right)
=I\left(  \beta\right).$$

 If
$j\in E,$ let $\left[  j\right]  _{\alpha}$ denote the $\sim_{\alpha}%
$-equivalence class of $j$, and let $E_{\alpha}$ denote the set of
all such equivalence classes. We can easily construct all of the
functions $\beta$ equivalent
to $\alpha$ in terms of injective functions%
\[
\sigma:E_{\alpha}\stackrel{1\mbox{-}1}{\rightarrow}H
\]
by defining
\[
\beta\left(  j\right)  =\sigma\left(  \left[  j\right]
_{\alpha}\right)  .
\]

Let $A$ contains exactly one function $\alpha$ from each
$\approx$-equivalence class of good functions. Then we can
write\begin{eqnarray}
I&=&\left|I\right|=\left|\frac{1}{k^2}\sum_{\alpha: E\rightarrow
H}\left( {\prod_{\nu=1}^{w}} \gamma_{\nu}\left(
\alpha({m_{\nu}\dot{+}1})\right) \right) \left(
{\prod_{\lambda=1}^{w}}\overline{ \gamma_{\lambda}\left(
\alpha((m_{\lambda}\dot{+}1)+m_w)\right) }
\right)I(\alpha)\right|\nonumber\\
&=& \left|\frac{1}{k^2}\sum_{\alpha\in A}I(\alpha)\sum_{\beta\approx
\alpha}\left( {\prod_{\nu=1}^{w}} \gamma_{\nu}\left(
\beta({m_{\nu}\dot{+}1})\right) \right) \left(
{\prod_{\lambda=1}^{w}}\overline{ \gamma_{\lambda}\left(
\beta((m_{\lambda}\dot{+}1)+m_w)\right) }
\right)\right|\nonumber\\
&=&\frac{1}{k^2}\left|\sum_{\alpha\in A}I(\alpha)\sum_{\sigma:
E_{\alpha}\rightarrow H}\left( {\prod_{\nu=1}^{w}}
\gamma_{\nu}\left( \sigma([m_{\nu}\dot{+}1]_{\alpha})\right) \right)
\left( {\prod_{\lambda=1}^{w}}\overline{ \gamma_{\lambda}\left(
\sigma([(m_{\lambda}\dot{+}1)+m_w]_{\alpha})\right) }
\right)\right|\nonumber\\
&\leq&\frac{1}{k^2}\sum_{\alpha\in A}\left|I\left(
\alpha\right)\right|\left|\sum_{\sigma:E_{\alpha
}\stackrel{1\mbox{-}1}{\rightarrow}H} \left( {\prod_{\nu=1}^{w}}
\gamma_{\nu}\left( \sigma([m_{\nu}\dot{+}1]_{\alpha})\right) \right)
\left( {\prod_{\lambda=1}^{w}}\overline{ \gamma_{\lambda}\left(
\sigma([(m_{\lambda}\dot{+}1)+m_w]_{\alpha})\right) }
\right)\right|.\end{eqnarray}

We know that \begin{equation}Card\left(  A\right)  \leq B\left(
2m\right) .\end{equation}

We only need to focus on $$ \left|\sum_{\sigma:E_{\alpha
}\stackrel{1\mbox{-}1}{\rightarrow}H} \left( {\prod_{\nu=1}^{w}}
\gamma_{\nu}\left( \sigma([m_{\nu}\dot{+}1]_{\alpha})\right) \right)
\left( {\prod_{\lambda=1}^{w}}\overline{ \gamma_{\lambda}\left(
\sigma([(m_{\lambda}\dot{+}1)+m_w]_{\alpha})\right) }
\right)\right|.
$$ Let
$$F_{\alpha}^1=\left\{\left[   m_ {\nu}\dot{+}1\right]
_{\alpha}:1\leq\nu\leq w,Card\left(\left[ m_ {\nu}\dot{+}1\right]
_{\alpha}\right) =1\right\},$$ $$F_{\alpha}^2= \left\{\left[(
 m_ {\nu}\dot{+}1)+m_w\right] _{\alpha}:1\leq\nu\leq
w,Card\left(\left[( m_ {\nu}\dot{+}1)+m_w\right] _{\alpha}\right)
=1\right\},$$ $$G_{\alpha}^1=\left\{ \left[  m_ {\nu}\dot{+}1\right]
_{\alpha}:1\leq \nu\leq w,Card\left( \left[  m_ {\nu}\dot{+}1\right]
_{\alpha}\right)
>1\right\},$$
$$G_{\alpha}^2=\left\{\left[(  m_ {\nu}\dot{+}1)+m_w\right] _{\alpha}:1\leq\nu\leq
w,Card\left(\left[( m_ {\nu}\dot{+}1)+m_w\right] _{\alpha}\right)
>1\right\},$$  $$K_{\alpha}=E_{\alpha}\backslash\left(
F_{\alpha}^1\cup F_{\alpha}^2\cup G_{\alpha}^1\cup G_{\alpha
}^2\right).$$ Since the product ${\prod_{\nu=1}^{w}}
\gamma_{\nu}\left( \sigma\left( \left[  m_ {\nu}\dot{+}1\right]
_{\alpha}\right) \right) $ is determined once $\sigma$ is defined on
$F_{\alpha}^1\cup F_{\alpha}^2\cup G_{\alpha}^1\cup G_{\alpha}^2,$
it follows that this product is repeated at most
$P\left(k,card\left( K_{\alpha}\right) \right) $ times. Hence we
have\begin{eqnarray}&& \left|\sum_{\sigma:E_{\alpha
}\stackrel{1\mbox{-}1}{\rightarrow}H} \left( {\prod_{\nu=1}^{w}}
\gamma_{\nu}\left( \sigma([m_{\nu}\dot{+}1]_{\alpha})\right) \right)
\left( {\prod_{\lambda=1}^{w}}\overline{ \gamma_{\lambda}\left(
\sigma([(m_{\lambda}\dot{+}1)+m_w]_{\alpha})\right) }
\right)\right| \nonumber\\
&\leq & P\left(k,card\left( K_{\alpha}\right) \right)\left| \sum
_{\sigma:F_{\alpha}^1\cup F_{\alpha}^2\cup G_{\alpha}^1\cup
G_{\alpha}^2\stackrel{1\mbox{-}1}{\rightarrow}F}\left(
{\prod_{\nu=1}^{w}} \gamma_{\nu}\left(
\sigma([m_{\nu}\dot{+}1]_{\alpha})\right) \right) \left(
{\prod_{\lambda=1}^{w}}\overline{ \gamma_{\lambda}\left(
\sigma([(m_{\lambda}\dot{+}1)+m_w]_{\alpha})\right) }
\right)\right| \nonumber\\
 &\leq&k^{card\left( K_{\alpha}\right)  }\left| \sum
_{\sigma:F_{\alpha}^1\cup F_{\alpha}^2\cup G_{\alpha}^1\cup
G_{\alpha}^2\stackrel{1\mbox{-}1}{\rightarrow}F}\left(
{\prod_{\nu=1}^{w}} \gamma_{\nu}\left(
\sigma([m_{\nu}\dot{+}1]_{\alpha})\right) \right) \left(
{\prod_{\lambda=1}^{w}}\overline{ \gamma_{\lambda}\left(
\sigma([(m_{\lambda}\dot{+}1)+m_w]_{\alpha})\right) }
\right)\right|.
\end{eqnarray}

If $a\in F_{\alpha}^1$ (or $F_{\alpha}^2$), from the definition of
$F_{\alpha}^1$ (or $F_{\alpha}^2$), the cardinality of $a$ is 1.
Then define $f_a\left( \sigma\left(
a\right)\right)=\gamma_{\nu}\left( \sigma\left( a\right)\right)$
(or$\overline{\gamma_{\nu}\left( \sigma\left( a\right)\right)}$). By
$\tau_k(x_i)=0$ for all $1\leq i\leq w$, it follows that
$\sum_{s=1}^k f_a(s)=0$. If $b\in G_{\alpha}^1$ (or $G_{\alpha}^2$),
from the definition of $G_{\alpha}^1$ (or $G_{\alpha}^2$), the
cardinality of $b$ is greater than 1, say it $r$. Then define
$g_b\left( \sigma\left(b\right)\right)=\left(\gamma_{\nu}
\left(\sigma\left( b\right)\right)\right)^r$ (or
$\overline{\left(\gamma_{\nu} \left(\sigma\left(
b\right)\right)\right)^r}$). Let
$F_{\alpha}=F_{\alpha}^1+F_{\alpha}^2$ and
$G_{\alpha}=G_{\alpha}^1+G_{\alpha}^2$. Then we have
\begin{eqnarray}&&\left| \sum
_{\sigma:F_{\alpha}\cup
G_{\alpha}\stackrel{1\mbox{-}1}{\rightarrow}F}\left(
{\prod_{\nu=1}^{w}} \gamma_{\nu}\left(
\sigma([m_{\nu}\dot{+}1]_{\alpha})\right) \right) \left(
{\prod_{\lambda=1}^{w}}\overline{ \gamma_{\lambda}\left(
\sigma([(m_{\lambda}\dot{+}1)+m_w]_{\alpha})\right) }
\right)\right|\nonumber\\
&=&\left|\sum _{\sigma:F_{\alpha}\cup
G_{\alpha}\stackrel{1\mbox{-}1}{\rightarrow}H}{\prod_{a\in
F_{\alpha}}}f_a\left(  \sigma\left( a\right)\right){\prod_{b\in
G_{\alpha}}}g_b\left( \sigma\left(
b\right)\right) \right|\nonumber\\
&&\left(\mbox{letting}\ F=F_{\alpha}, G=G_{\alpha}\ \mbox{and using Lemma \ref{0sums} }\right)\nonumber\\
&\leq & k^{\left[  card\left( F_{\alpha}\right)  /2\right]
+card\left( G_{\alpha}\right) }(2w)^{2w}M^{2w}.
\end{eqnarray}
As we mentioned before that $card\left( \left[ j\right]
_{\alpha}\right) =1$ implies $\left[  j\right] _{\alpha}\in
F_{\alpha},$ we see that
\begin{equation} \left[  card\left( F_{\alpha}\right) /2\right] +card\left(
G_{\alpha }\right) +card\left( K_{\alpha}\right)  \leq card\left(
E\right) /2=2m_{w}/2=m.
\end{equation}
Combining (9), (10), (11), (12), (13) and (14) together, we have $$
|I|\leq \frac{1}{k^2}B\left(  2m\right)  4^{m^{2}}\left( 2Mw\right)
^{2w}.$$ This completes the proof of the second statement.

{\it Proof of the third statement.}  The third statement follows
from statement 1 and statement 2 and Chebychev's inequality.

Let $A=B\left( m\right)  2^{m^{2}}\left( Mw\right) ^{w}$ and
$B=B\left(  2m\right)  4^{m^{2}}\left( 2Mw\right) ^{2w}$. Define
$$f(\overrightarrow{v})=\tau_{k}\left(  g_{1}\left( \vec {u}\right)
x_{1}g_{2}\left(  \vec{u}\right)  x_{2}\cdots g_{w}\left(  \vec
{u}\right)  x_{w}\right). $$ Then the expected value of $f$
$E(f)=\int_{\mathcal{U}_{k}^{n}}f(\overrightarrow{v})d\mu_k^n(\overrightarrow{v}),$
and the variance of $f$ is
\begin{eqnarray*}\mbox{Var}(f)&=&\int_{\mathcal{U}_{k}^{n}}\left\vert
\tau_{k}\left( g_{1}\left( \vec {v}\right)  x_{1}g_{2}\left(
\vec{v}\right) x_{2}\cdots g_{w}\left(  \vec {v}\right) x_{w}\right)
\right\vert ^{2}d\mu_{k}^{n}\left( \vec{v}\right)\\
&&- \left\vert\int_{\mathcal{U}_{k}^{n}}\tau_{k}\left( g_{1}\left(
\vec {v}\right) x_{1}g_{2}\left(  \vec{v}\right) x_{2}\cdots
g_{w}\left( \vec {v}\right)  x_{w}\right) d\mu_{k}^{n}\left(
\vec{v}\right)
\right\vert^2\\
&\leq &\int_{\mathcal{U}_{k}^{n}}\left\vert \tau_{k}\left(
g_{1}\left( \vec {v}\right)  x_{1}g_{2}\left( \vec{v}\right)
x_{2}\cdots g_{w}\left(  \vec {v}\right) x_{w}\right)  \right\vert
^{2}d\mu_{k}^{n}\left( \vec{v}\right)\\
&\leq & \frac{B}{k^2}.\end{eqnarray*} Since
\begin{eqnarray*}\left|f(\overrightarrow{v}-E(f))\right|\geq
\left|f(\overrightarrow{v})\right|-\left|E(f)\right|\geq
\varepsilon-\frac{A}{k}>\frac{\varepsilon}{2},\end{eqnarray*} we
have $\{\overrightarrow{v}: |f(\overrightarrow{v})|\geq
\varepsilon\}\subseteq \{\overrightarrow{v}:
|f(\overrightarrow{v})-E(f)|\geq \frac{\varepsilon}{2}\}.$ Therefore
form Chebychev's inequality, we have
$$\mu_{k}^n(\{\overrightarrow{v}: |f(\overrightarrow{v})|\geq
\varepsilon\})\leq \frac{Var(f)}{\varepsilon^2}\leq
\frac{4B}{k^2\varepsilon^2}.$$  \hfill$\Box$

The following corollary if a direct consequence of the third
statement of Theorem \ref{theorem,estimation}.
\begin{corollary} Suppose $M, N, k$ are positive integers. Let $\cal D$ be a finite set of
commuting normal matrices with trace $0$
  in ${\cal M}_k(\Bbb C)$ and $\|x\|\leq M$ for all $x\in \cal D$.
  Let \begin{eqnarray*}{\cal E}&=&\{(g_1, \ldots, g_r, x_1, \ldots, x_r): r\in {\Bbb N}, g_1, \ldots, g_r\ \mbox{are reduced words in}\ {\Bbb F}_n\setminus \{e\}\\
&& \mbox{such that}\
 \sum_{i=1}^r\mbox{length}(g_i)\leq N,\ \mbox{and}\ x_1, \ldots , x_r \in {\cal
 D}\}.\end{eqnarray*}
 Suppose $\mathfrak{e}=(g_1, \ldots, g_r, x_1, \ldots, x_r)\in \cal E$,
 define $\mathfrak{e}(\overrightarrow{v})=g_1(\overrightarrow{v})x_1\cdots
 g_r(\overrightarrow{v})x_r$. Then $$\mu_k^n\left(\bigcap_{\mathfrak{e}\in {\cal E}}
 \left\{\overrightarrow{v}: \left|\tau_k(\mathfrak{e}(\overrightarrow{v}))\right|<\varepsilon\right\}\right)\geq\ 1-card({\cal E})\frac{4B\left(  2m\right)
   4^{m^{2}}\left( 2Mw\right) ^{2w}}{k^2\varepsilon^2}.$$
\end{corollary}

\vspace{0.2cm}

Lemma 5.1 \cite{Vo1} follows directly from the corollary above.

Let $\cal M$ be a von Neumann algebra with a tracial state $\tau$
and $X_1, X_2, \ldots, X_n$ be elements in ${\cal M}$. For any
$R,\varepsilon>0$, and positive integers $m$ and $k$, define
$\Gamma_R(X_1,\ldots,X_n;m,k,\varepsilon)$ to be the subset of
${\cal M}_k({\Bbb C})^n$ consisting of all $(x_1,\ldots,x_n)$ in
${\cal M}_k({\Bbb C})^n$ such that $\|x_j\|\leq R$ for $1\leq j\leq
n$, and
$$|\tau_k(x_{i_1}^{\eta_1}\cdots
x_{i_q}^{\eta_q})-\tau(X_{i_1}^{\eta_1}\cdots
X_{i_q}^{\eta_q})|<\varepsilon,$$ for all $1\leq i_1,\ldots,i_q\leq
n$, all $\eta_1,\ldots,\eta_q\in\{1,*\}$ and all $q$ with $1\leq
q\leq m$.

Suppose $\overrightarrow{U}$ is a n-tuple in $\cal M$ and, for each
positive integer $k$, $\overrightarrow{u_k}$ is a n-tuple in ${\cal
M}_k(\Bbb C)$, then we say {\it $\overrightarrow{u_k}$ converges to
$\overrightarrow{U}$ in distribution} if
$p(\overrightarrow{u_k})\rightarrow p(\overrightarrow{U})$ for all
$*$-monomials $p.$

\begin{corollary} Let $M,m$ be positive integers and $\varepsilon>0$. Suppose $\cal M$ is a von Neumann algebra with a faithful trace $\tau$.
Suppose $X_1, \ldots, X_s$ are commuting normal operators in $\cal
M$, $U_1, \ldots, U_n$ are free Haar unitary elements in $\cal M$
and $\{X_1, \ldots, X_s\},\{U_1, \ldots, U_n\}$ are free. For any
positive integer $k$, let $\{x(k,1), \ldots, x(k,s)\}$ be a set of
commuting normal $k\times k$ matrices such that
$\sup_{k,j}\|x(k,j)\|\leq M$ and
$$(x(k,1), \ldots, x(k,s))\rightarrow (X_1, \ldots, X_s)$$ in distribution as $k\rightarrow \infty$.

 If
$$\Omega_k=\{(v_1, \ldots, v_n)\in {\cal U}_k^n: (x(k,1), \ldots, x(k,s),v_1, \ldots, v_n)\in\Gamma_1(X_1, \ldots, X_s, U_1, \ldots, U_n; m, k,\varepsilon)\},$$
then $$\lim_{k\rightarrow \infty}\mu_k^n(\Omega_k)=1.$$
\end{corollary}

 Lemma 5.2 \cite{Vo1} follows directly from the corollary above.

We end this paper with one last corollary.
\begin{corollary}Let $M,m$ be positive integers and $\varepsilon>0$. Suppose $\cal M$ is a von Neumann algebra with a faithful trace $\tau$.
Suppose $X_1, \ldots, X_s$ are free normal operators in $\cal M$.
Suppose $\{x(k,1), \ldots, x(k,s)\}$ is a set of normal $k\times k$
matrices such that $\sup_{k,j}\|u(k,j)\|\leq M$ and, for $1\leq
j\leq s$, $x(k,j)\rightarrow X_j$ in distribution as $k\rightarrow
\infty$.

 If
$$\Theta_k=\{(v_1, \ldots, v_s)\in {\cal U}_k^s: (v_1^*x(k,1)v_1, \ldots, v_s^*x(k,s)v_s)\in\Gamma_1(X_1, \ldots, X_s; m, k,\varepsilon)\},$$
then $$\lim_{k\rightarrow \infty}\mu_k^n(\Theta_k)=1.$$
\end{corollary}

\end{document}